\documentclass[12pt]{article}
\usepackage{amsmath, amssymb,amscd, amsthm}

\newtheorem{theorem}{Theorem}[section]

\newtheorem{cor}[theorem]{Corollary} 
\newtheorem{prop}[theorem]{Proposition} 

\theoremstyle{definition}

\theoremstyle{remark}
\newtheorem{remark}[theorem]{Remark}
\newtheorem{rems}[theorem]{Remarks} 

\numberwithin{equation}{section}

\begin{document} 

\title 
{Sigma function solution of the initial value problem  
for Somos 5  
sequences} 

\author{
A.N.W. Hone\thanks{   
Institute of Mathematics, Statistics \&
Actuarial Science,  
University of Kent,  
Canterbury CT2 7NF, UK.    
email:  anwh@kent.ac.uk} }  



\newcommand{\beq}{\begin{equation}}
\newcommand{\eeq}{\end{equation}}
\newcommand{\bea}{\begin{eqnarray}}
\newcommand{\eea}{\end{eqnarray}}
\newcommand\la{{\lambda}} 
\newcommand\ka{{\kappa}} 
\newcommand\al{{\alpha}} 
\newcommand\be{{\beta}} 
\newcommand\om{{\omega}}
\newcommand\tal{{\tilde{\alpha}}}
\newcommand\tbe{{\tilde{\beta}}} 
\newcommand\tla{{\tilde{\lambda}}}
\newcommand\tmu{{\tilde{\mu}}}
\newcommand\si{{\sigma}}
\newcommand\lax{{\bf L}}  
\newcommand\mma{{\bf M}}  
\newcommand\rd{{\mathrm{d}}}
\newcommand\tI{{\tilde{I}}}
\newcommand\tJ{{\tilde{J}}}


\maketitle 

\begin{abstract}
The Somos 5 sequences are 
a family of sequences defined by a fifth order 
bilinear recurrence relation with constant coefficients. For 
particular choices of coefficients and initial data, integer 
sequences arise. By making the connection with a 
second order nonlinear mapping with a first integral, we 
prove that the two subsequences of odd/even index terms
each satisfy a Somos 4 (fourth order) recurrence. This leads 
directly to  
the explicit solution of the initial value problem 
for the Somos 5 sequences in terms of the Weierstrass sigma 
function for an associated elliptic curve. 
\end{abstract}

\maketitle
 
\section{Introduction}

The subject of this article is the bilinear 
recurrence relation 
\beq 
\tau_{n+3}\tau_{n-2}=\tal \,\tau_{n+2}\tau_{n-1}+\tbe \,\tau_{n+1}\tau_n,  
\label{s5} 
\eeq 
where $\tal$ and $\tbe$ are constant parameters, which is 
known as the Somos 5 recurrence \cite{gale, rob, sloane}. Somos 5 sequences 
are defined by specifying initial data $\tau_0$, 
$\tau_1$, $\tau_2$, $\tau_3$, $\tau_4$ together with the values of the 
coefficients $\tal$, $\tbe$, and when all of these are
suitably chosen 
then integer sequences can result. For example, 
when 
\beq \label{s5data} 
\tau_0=\tau_1=\tau_2=\tau_3=\tau_4=1, \qquad \tal=\tbe =1 
\eeq 
then the integer sequence  
\beq 
\label{somos5} 
1,1,1,1,1,2,3,5,11,37,83,274,1217,6161,22833,\ldots 
\eeq 
results. Traditionally the latter is referred to as {\it the} Somos 5 
sequence, but we will follow the terminology of \cite{rob, swart}  
and refer to this specific sequence as Somos (5).  

In recent work  \cite{honeblms}, we 
have considered fourth
order quadratic recurrences of the form
\beq
\tau_{n+2}
\tau_{n-2}
=\al \,
\tau_{n+1}
\tau_{n-1}
+\beta
\,
(\tau_n)^2,
\label{bil}
\eeq
where $\al$ and $\beta$ are constant parameters. Sequences 
that satisfy such 
fourth order recurrences are known as Somos 4 sequences
\cite{rob, swart}, and they originally arose  
in the theory of elliptic divisibility sequences
\cite{ward1, ward2, shipsey}. 
In that context, both the parameters $\al$, $\beta$
and the iterates $\tau_n$
are usually integers, 
in which case the sequences have applications in
number theory, as they provide a potential source
of large prime numbers \cite{eew, ems}. Moreover Somos 4 
sequences  
provide a simple example  
of the Laurent phenomenon: taking the initial data
$\tau_0,\tau_1,\tau_2,\tau_3$ and the parameters $\al ,\beta$ as
variables, all subsequent terms $\tau_n$ for $n\geq 4$ in the sequence
are Laurent polynomials 
(rather than just rational functions) in these variables.
Fomin and Zelevinsky have proved
that this remarkable ``Laurentness'' property is shared by a variety of
other recurrences in one and more dimensions, including 
the Somos 5 recurrences (\ref{s5}), and this has applications in 
combinatorics and commutative algebra (see \cite{fz} and references). 
Although Somos sequences are often defined over $\mathbb{Q}$, 
here we will work over $\mathbb{C}$ and  
occasionally indicate when  our results restrict  
to yield rationals or integers.

Morgan Ward first introduced elliptic divisibility 
sequences \cite{ward1, ward2},  
a family of antisymmetric sequences with $a_n=-a_{-n}$ 
defined by recurrences
of the form
\beq
\label{edseq}
a_{n+2}a_{n-2}=(a_2)^2a_{n+1}a_{n-1}
-a_1 a_3(a_n)^2,
\eeq
by considering sequences of rational points $nP$
on an elliptic curve $E$ over $\mathbb{Q}$ 
(see
also \cite{recs, shipsey}).  
To obtain integer
sequences of this kind it is required that
\beq \label{intres} 
a_0=0, \quad a_1=1, \quad a_2,a_3,a_4\in\mathbb{Z}
\qquad \mathrm{with}\quad a_2|a_4,
\eeq 
and then it turns out that the subsequent terms of the sequence 
satisfy the remarkable divisibility property
\beq
\label{divis}
a_n  |a_m \qquad \mathrm{ whenever}
\qquad n|m.
\eeq
These integrality and divisibility property can be said to 
rest on 
the fact that the terms of an elliptic divisibility 
sequence satisfy the 
Hankel determinant relation
\beq \label{hankel}
a_{n+m}\, a_{n-m}=\left|\begin{array}{cc} a_m a_{n-1} &
a_{m-1}a_n \\
a_{m+1}a_{n} &
a_{m}a_{n+1} \end{array}\right|,
\eeq
for all $m,n\in\mathbb{Z}$. 
Starting from the Hankel determinant formula, with initial 
data of the form (\ref{intres}), it is easy to
prove by induction that all $a_n$ for $n\geq 0$ are integers with the
divisibility property (\ref{divis}). Although with $a_0=0$ 
there is potentially the problem 
of dividing by zero, it turns out that these sequences are 
consistently extended to all $n\in\mathbb{Z}$.  

As we explain below, 
the formula (\ref{hankel}) itself is most easily proved by making 
use of Ward's 
explicit expression 
\beq
\label{edsform}
a_n=\frac{\si (n\ka )}{\si (\ka )^{n^2}}.
\eeq
for the general term of the sequence in
terms of the Weierstrass sigma function associated with the curve $E$, 
and then applying the addition formula 
\beq
\frac{\si (z+\ka )\,\si (z-\ka )}
{\si (z)^2\si (\ka )^2}
=\wp (\ka )-\wp (z )
\label{addweier}
\eeq
or its immediate consequence, 
\bea  \label{3term} \nonumber  
\sigma (c+d) \,\sigma (c-d) \,\sigma (a+b) \,\sigma (a-b) && \\  
- 
\sigma (b+d) \,\sigma (b-d) \,\sigma (a+c) \,\sigma (a-c) && \\ 
+ 
\sigma (b+c) \,\sigma (b-c) \,\sigma (a+d) \,\sigma (a-d) & = & 0, \nonumber  
\eea 
which is called the three-term equation  (see 
$\S$20.53 in \cite{ww}).  
Ward's elliptic divisibility sequences are just special 
cases of Somos 4 sequences. 
The following theorem was proved in 
\cite{honeblms}. 


\begin{theorem} 
The general solution
of the Somos 4 recurrence relation
(\ref{bil}) 
takes the form
\beq
\tau_n=
A\,B^n\frac{\si (z_0 +n\ka )}{
\si (\ka )^{n^2}},
\label{form}
\eeq
where $\ka$ and $z_0$ are non-zero
complex numbers, the constants $A$ and $B$ are given by
\beq
A=\frac{\tau_0}{\si (z_0 )}, \qquad
B=\frac{\si (\ka )\si (z_0 )\,\tau_1}
{\si (z_0+\ka )\,\tau_0},
\label{consts}
\eeq
and $\si (z)=\si (z;g_2,g_3)$ denotes the Weierstrass sigma function
of an associated elliptic curve
\beq \label{canon}
E: \qquad y^2=4x^3-g_2x-g_3.
\eeq
The values $\ka$, $z_0$  and the invariants $g_2$, $g_3$ are precisely
determined from the initial data
$\tau_0,\tau_1,\tau_2,\tau_3$ and the parameters $\al ,\beta$, 
where the latter are given as elliptic functions of $\ka$ by 
\beq \label{abform} 
\al =\wp '(\ka )^2, \qquad 
\be = \wp '(\ka )^2 \Big(\wp (2\ka )-\wp (\ka )\Big) . \eeq 
\end{theorem}

The terms of a Somos 4 sequence correspond to 
the sequence of points 
$P_0+nP$ on the curve $E$, with $z_0+n\ka $ being the 
associated sequence of points on the Jacobian. 
The elliptic divisibility 
sequences arise when $P_0\to\infty$ (or 
equivalently $z_0\to 0$ with $\tau_0\to 0$). 
To illustrate the effectiveness of such explicit formulae for these 
sequences  
in terms of Weierstrass functions, we can mention two corollaries 
of the above theorem. 


\begin{cor} 
The terms of a Somos 4 sequence satisfy the 
Hankel determinant formula   
\beq 
\label{hankel2} 
\tau_{n+m}\tau_{n-m}= 
\left|\begin{array}{cc} 
a_m \tau_{n+1} & a_{m-1} \tau_n \\ 
a_{m+1}\tau_n & a_m \tau_{n-1} 
\end{array} \right|  
\eeq 
for all $m,n\in \mathbb{Z}$, where $a_m=\si (m\ka )/\si 
(\ka )^{m^2}$ are the terms of an associated elliptic divisibility 
sequence.  
\end{cor} 


\begin{proof}  
Substituting the explicit expressions (\ref{edsform}) and (\ref{form}) 
for $a_n$ and $\tau_n$ into (\ref{hankel2}), setting $z=z_0+n\ka$ 
and using (\ref{addweier}) yields 
$$ 
\frac{\si (\ka )^{2m^2+2n^2}}{A^2B^{2n}\si (m\ka )^2\si (z)^2} 
\Big\{
\tau_{n+m}\tau_{n-m}-(a_m)^2\tau_{n+1}\tau_{n-1}+a_{m+1}a_{m-1}(\tau_n)^2 
\Big\} $$ 
$$=\frac{\si (z+m\ka )\si (z-m\ka )}{\si (z)^2\si (m\ka )^2}   
-\frac{\si (z+\ka )\si (z-\ka )}{\si (z)^2\si (\ka )^2}  
+\frac{\si ((m+1)\ka )\si ((m-1)\ka )}{\si (\ka )^2\si (m\ka )^2} 
$$ 
$$ 
=\wp (m\ka ) - \wp (z) -\wp (\ka )+\wp(z)+\wp(\ka )-\wp(m\ka ) 
 =0. 
$$  
\end{proof} 



\begin{cor}    
The terms 
of a Somos 4 sequence also satisfy 
\beq 
\label{hankel3}
a_1 a_2 \tau_{n+m+1}\tau_{n-m}=
\left|\begin{array}{cc}
a_{m+1} \tau_{n+2} & a_{m-1} \tau_n \\
a_{m+2}\tau_{n+1} & a_m \tau_{n-1}
\end{array} \right| 
\eeq  
for all $m,n\in \mathbb{Z}$, with $a_n$ for 
$n\in \mathbb{Z}$ given by (\ref{edsform})  
as above.    

\end{cor} 


\begin{proof}  

Substituting the formulae (\ref{edsform}) and (\ref{form})
for $a_n$ and $\tau_n$ into (\ref{hankel3}), and setting 
$\hat{z}=z_0+(n+\frac{1}{2})\ka $  
gives 
$$ 
\frac{\si (\ka )^{2m^2+2n^2+2m+2n+6}}{A^2B^{2n+1}}
\Big\{
a_1a_2\tau_{n+m+1}\tau_{n-m} 
-a_m a_{m+1}\tau_{n+2}\tau_{n-1}+a_{m-1}a_{m+2}\tau_{n+1}\tau_n
\Big\} $$
$$ \begin{array}{lcc} 
= \,\si (\ka ) \,\si (2\ka ) \,\si \Big(\hat{z}+(m+1/2)\ka \Big) 
\,\si \Big(\hat{z}-(m+1/2)\ka \Big) && \\  
-\,\si (m\ka ) \,\si \Big((m+1)\ka \Big) \,\si \Big(\hat{z}+3\ka /2 \Big)
\,\si \Big(\hat{z}-3\ka /2 \Big) && \\  
+\,\si \Big((m-1)\ka \Big) \,\si \Big((m+2)\ka \Big) 
\,\si \Big(\hat{z}+\ka /2\Big)
\,\si \Big(\hat{z}-\ka /2\Big)
&= &0  
\end{array}   
$$  
by (\ref{3term}).  

\end{proof} 



\begin{rems}  
These two corollaries of Theorem 1.1 
are the main result of \cite{vdpswart}, 
where they are proved by purely algebraic  
means i.e. without using explicit analytic 
expressions such as (\ref{form}). 
Morgan Ward's formula (\ref{hankel}) for elliptic divisibility 
sequences is clearly just a special case of the first corollary. Comparing 
the statement of Corollary 1.2 when $m=2$ with the expressions for 
the coefficients of the Somos 4 recurrence in (\ref{abform}) gives 
$\al=(a_2)^2$, $\be =-a_1a_3$, so it is clearly 
required that 
\beq \label{efns} 
\wp '(\ka )^2=\frac{\si (2\ka )^2}{\si (\ka )^8}, 
\qquad 
\wp '(\ka )^2\left( \wp (2\ka )-\wp (\ka )\right) 
=-\frac{\si (3\ka )}{\si (\ka ) 
^9}.    
\eeq 
To see why these two identities must hold, note that for any $n$,
$a_n=\si (n\ka )/\si (\ka )^{n^2}$ is an elliptic function of $\ka$ 
(see exercise 24 in Chapter 20 of \cite{ww}). 
Thus in the first formula of (\ref{efns}), both sides are elliptic 
functions of $\ka$ whose Laurent expansions around $\ka =0$ have the same 
principal parts; moreover both functions vanish at any one of the 
half-periods $\ka =\om_j$, $j=1,2,3$; 
hence these two expressions for $\al$ are
equivalent, and thence the two different formulae for $\be$ are seen to be 
equivalent by taking $z=2\ka$ in (\ref{addweier}).  
\end{rems} 

The main results of this paper are found in the next section: 
Theorem 2.7 presents the explicit solution of a second order nonlinear 
mapping associated with Somos 5 sequences, while 
Theorem 2.9  
presents the explicit form of the general solution for these sequences.  
In the third section we briefly illustrate the solution of the  
initial value problem for such sequences, by applying  
the theorem to the specific example of the Somos (5)  
sequence (\ref{somos5}), 
which is sequence A006721 in Sloane's catalogue \cite{sloane}.

\section{General solution of Somos 5}

The Somos 4 recurrence has a simple gauge invariance 
property: if $\tau_n$ satisfies (\ref{bil}) then for any 
non-zero constants $\tilde{A},\tilde{B}$ it is clear that 
the transformation  
\beq\label{gauge} 
\tau_n\longrightarrow\tilde{\tau}_n= \tilde{A}\, \tilde{B}^n\,\tau_n 
, \qquad n\in\mathbb{Z} 
\eeq 
yields  another solution of  
the same fourth order recurrence. Thus it is natural to introduce 
the gauge-invariant combination 
\beq 
f_n=\frac{\tau_{n+1}\tau_{n-1}}{(\tau_n)^2}, 
\label{fdef} 
\eeq 
which satisfies the second order nonlinear mapping 
\beq 
f_{n+1}=\frac{1}{f_{n-1}f_n}\left( \al +\frac{\be}{f_n}\right) 
\label{nonlia} 
\eeq   
(written in the more symmetrical form (\ref{nonli}) below).
The following result in \cite{honeblms} 
was the basis for 
the proof of Theorem 1.1. 
  

\begin{prop} 
The second order nonlinear mapping 
\beq 
f_{n-1}(f_n)^2f_{n+1}=\al f_n +\be 
\label{nonli} 
\eeq 
has the conserved quantity 
\beq 
J=f_{n-1}f_n +\al \left( \frac{1}{f_{n-1}}+\frac{1}{f_n}
\right)+ \frac{\be}{f_{n-1}f_n}. 
\label{jint} 
\eeq 
Its general solution has the form  
\beq 
f_n=\wp (\ka )-\wp (z_0 +n\ka ),  
\label{desol} 
\eeq 
where $\wp (z)=\wp(z;g_2,g_3)$ is the  Weierstrass elliptic function 
associated with an elliptic curve (\ref{canon}) 
with invariants $g_2$, $g_3$. 
The parameters $z_0$, 
$\ka\in\mathbb{C}$ and the two invariants are determined from 
$\al$, $\be$ and the initial data $f_0$, $f_1$ 
according to the formulae 
\beq \label{lamf} 
\la =\frac{1}{3\al}\left( \frac{J^2}{4}-\be \right), \quad 
g_2=12\la^2-2J, \quad g_3=4\la^3-g_2\la -\al , 
\eeq 
together with the elliptic integrals 
\beq 
z_0=\pm\int_\infty^{\la -f_0} \frac{\rd x}{y}, 
\qquad 
\ka =\pm\int_\infty^{\la } \frac{\rd x}{y}. 
\label{eints} 
\eeq
\end{prop}

\begin{rems} 
The above is a slight 
reformulation of Proposition 2.2 in \cite{honeblms}, 
where 
the quantity $\la =\wp (\ka )$ was taken as the basic 
conserved quantity 
rather than the simpler quantity 
$J=\wp ''(\ka )$; the latter is the same as the ``translation 
invariant'' found in Swart's thesis \cite{swart}.  
There is an ambiguity 
in the sign of the integrals (\ref{eints}) due to the 
elliptic involution $y\to -y$, but the solution (\ref{desol}) 
is clearly invariant 
under $z_0\to -z_0$, $\ka\to -\ka$ due to the fact that $\wp(z)$ is 
an even function; the relative sign of $z_0$ and 
$\ka$ is fixed by the consistency requirement 
$\wp '(z_0)\wp' (\ka)=(f_0)^2(f_{1}-f_{-1})$. 
The mapping (\ref{nonli}), or equivalently 
(\ref{nonlia}), is of the type studied 
by Quispel, Roberts and Thompson \cite{qrt} (see also 
\cite{bastien, iatrou1, iatrou2, tsuda} and references). 
The solution (\ref{desol}) of (\ref{nonli})  
corresponds to the solution (\ref{form}) of the 
Somos 4 recurrence,  
modulo the two-parameter abelian group of gauge 
transformations (\ref{gauge}). 
\end{rems} 


\begin{prop} 
By means of the substitution (\ref{fdef}), 
the Somos 5 recurrence (\ref{s5}) corresponds to 
the third order nonlinear mapping 
\beq 
f_{n-1}(f_n)^2(f_{n+1})^2f_{n+2}=\tal f_nf_{n+1} +\tbe ,  
\label{nonli3}  
\eeq 
which has two independent conserved quantities 
$\tI$, $\tJ$ given by  
\beq 
\tI=f_{n-1}f_nf_{n+1} +\tal \left( \frac{1}{f_{n-1}}+\frac{1}{f_n} 
+\frac{1}{f_{n+1}}\right)+ \frac{\tbe}{f_{n-1}f_nf_{n+1}}, 
\label{i3int} 
\eeq 
\beq 
\tJ=f_{n-1}f_n+f_nf_{n+1} +\tal 
\left( \frac{1}{f_{n-1}f_n}+\frac{1}{f_nf_{n+1}}
\right)+ \frac{\tbe}{f_{n-1}(f_n)^2f_{n+1}}. 
\label{j3int} 
\eeq 
Any solution of the Somos 4 recurrence (\ref{bil}) 
also satisfies a Somos 5 recurrence (\ref{s5}) with parameters 
\beq \label{s4ins5} 
\tal =-\be , \qquad \tbe =\al^2 +\be J, 
\eeq 
where $J$ is as in (\ref{jint}), and values of the 
conserved quantities 
for the associated 
solution of (\ref{nonli3}) are   
\beq \label{ijvals} 
\tJ =J, \qquad \tI =2\al. 
\eeq
\end{prop} 


\begin{proof} 
Letting $\tI_n$ denote the expression on the right hand side of 
(\ref{i3int}), and letting $\tI_{n+1}$ be the same expression 
with $n\to n+1$, we have immediately that $\tI_{n+1}-\tI_n =0$ whenever 
(\ref{nonli3}) holds, and similarly $\tJ$ in (\ref{j3int}) is conserved 
by the third order nonlinear mapping. Now for any sequence 
$\tau_n$, $n\in\mathbb{Z}$ that is a 
solution of a 
Somos 4 recurrence (\ref{bil}), there is a corresponding solution 
of (\ref{nonli}) given by (\ref{fdef}). By considering the expression 
$$ 
\mathcal{E}_n:= f_{n-1}(f_n)^2(f_{n+1})^2f_{n+2}-\tal f_nf_{n+1} -\tbe , 
$$ 
and then using (\ref{nonli}), first  with  $n\to n+1$ to eliminate $f_{n+2}$, 
and then once more to eliminate $f_{n+1}$, we find 
$$ 
\mathcal{E}_n=\al^2-\tbe +\be f_nf_{n-1}+\be \left( 
\frac{\al}{f_n}-\frac{\tal}{f_nf_{n-1}}\right) 
-\frac{\al\tal}{f_{n-1}} 
$$ 
$$ 
=\al^2 +\be J-\tbe  
-(\be +\tal )\left( \frac{\al}{f_{n-1}}+\frac{\be}{f_{n-1}f_n
}  
\right) 
$$ 
(where we have used (\ref{jint}) at the last step) 
and hence $\mathcal{E}_n=0$ for all $n$ only when the 
two conditions (\ref{s4ins5})  
hold. In that case, $f_n$ satisfies (\ref{nonli3}) with these values of    
$\tal$, $\tbe$, which is equivalent to  $\tau_n$ being a solution of 
(\ref{s5}). The values (\ref{ijvals}) follow from similar calculations. 
\end{proof}

\begin{remark} 
Clearly the special case $m=2$ of 
Corollary 1.2 above means that every Somos 4 is a Somos 5 with 
$\tal =  a_3$, $\tbe = -a_4/a_2$, so using (\ref{edsform}) 
to compare with (\ref{abform}) and noting that $J=\wp '' (\ka )$ 
it is clear that the equation $\tal = a_3 =-\be$ is the same as 
the first equation of (\ref{efns}), while the formula for 
$\tbe$ in (\ref{s4ins5}) yields the interesting identity 
\beq \label{interest} 
\wp '(\ka )^4 +\wp '' (\ka ) \wp '(\ka )^2\Big(\wp (2\ka ) 
 -\wp (\ka )\Big) =-\frac{\si (4\ka )}{\si (2\ka ) \si (\ka )^{12}}. 
\eeq 
The latter should be very tedious to prove by comparing poles and zeros, 
since on each side the leading term in the   
Laurent expansion at $\ka =0$   is $-2/\ka^{12}$. 
\end{remark} 

The formula (\ref{form}) is obviously not general enough to capture  
all of the Somos 5 
sequences: not every Somos 5 is a Somos 4! 
The expression (\ref{form}) only depends on the six parameters 
$A$, $B$, $z_0$,
$\ka$, $g_2$, $g_3$, whereas to specify a generic Somos 5 sequence 
requires seven parameters, corresponding to the five initial data  
$\tau_0$, $\tau_1$, $\tau_2$, $\tau_3$, $\tau_4$ plus the two  
coefficients $\tal$, $\tbe$. One way to interpret the extra 
degree of freedom 
in Somos 5 sequences is by noting that they have a three-parameter
group of gauge transformations given by  
\beq\label{s5gauge} 
\tau_{2n}\longrightarrow\tilde{\tau}_{2n}= 
\tilde{A}_+\, \tilde{B}^{2n}\,\tau_{2n},  
\quad  
\tau_{2n+1}\longrightarrow\tilde{\tau}_{2n+1}= 
\tilde{A}_-\, \tilde{B}^{2n+1}\, 
\tau_{2n+1}, \quad n\in\mathbb{Z}.  
\eeq  
In other words, compared with Somos 4 
there is an extra freedom to rescale the even/odd terms differently while 
leaving the recurrence (\ref{s5}) unchanged. The quantity $f_n$ 
defined by (\ref{fdef}) is not invariant under this larger group 
of transformations, and thus it is natural to consider something that is, 
namely  
\beq \label{hdef}  
h_n=f_{n+1}f_n=\frac{\tau_{n+2}\tau_{n-1}}{\tau_{n+1}\tau_n}.   
\eeq   

Before we consider the recurrence satisfied by $h_n$, it 
is convenient to state the following general result on 
mappings of the plane associated with biquadratic curves, which have been 
studied in detail recently by Iatrou and Roberts \cite{iatrou1, iatrou2},  
Bastien and Rogalski \cite{bastien} and also Tsuda \cite{tsuda}.   

\begin{prop} 
The quartic curve $\hat{E}$ defined by  
\beq \label{biq} 
\mathcal{B}(X,Y)\equiv 
eX^2Y^2+dXY(X+Y)+c(X^2+Y^2)+b(X+Y)+a-KXY=0 
\eeq 
has genus one. The mapping of the plane 
\beq \label{plamap}  
(u_{n-1},u_n)\mapsto (u_n,u_{n+1}) 
\eeq 
defined by 
\beq \label{bimap} 
u_{n+1}u_{n-1}= 
\frac{a+bu_n+cu_n^2}{c+du_n+eu_n^2} 
\eeq 
corresponds to the composition of two involutions on $\hat{E}$, 
and has the conserved quantity 
\beq 
\tilde{K}=eu_{n-1}u_n+d(u_{n-1}+u_n)+c\left( 
\frac{u_{n-1}}{u_n}+\frac{u_{n}}{u_{n-1}}\right) 
+b\left(\frac{1}{u_{n-1}}+\frac{1}{u_n}\right)+ 
\frac{a}{u_{n-1}u_n}. 
\label{kint} 
\eeq 
The iterates of (\ref{bimap}) correspond  
to the the sequence of points 
\beq \label{pseq}  
\hat{P}_n = (u_0,u_1) \, + \, n\, \hat{P}\,\in \, \hat{E}  
\eeq 
for some $\hat{P}\in\hat{E}$. Furthermore, upon 
uniformizing $\hat{E}$,  
the 
coordinates $u_n$ can be written as  
\beq 
u_n  = f(z_0+n\ka ) \label{efn}  
\eeq 
in terms of an elliptic function $f$ and certain parameters 
$z_0,\ka \in\mathbb{C}$. 
\end{prop} 

\begin{proof} 
To see that the curve $\hat{E}$ is elliptic, it is 
sufficient to observe that by the transformation 
$\mathcal{R}=XY$, $\mathcal{S}=X+Y$, this symmetric 
biquadratic is a ramified double cover  of the rational 
curve  
$$ 
e\mathcal{R}^2+d\mathcal{R}\mathcal{S}+c(\mathcal{S}^2-2\mathcal{R}) 
+b\mathcal{S}+a-K\mathcal{R}=0 
$$ 
with branching at four points $(X,X)$ given by the 
roots of the quartic $\mathcal{B}(X,X)=0$, and so by 
the Hurwitz formula the genus is one.  
The symmetric biquadratic
(\ref{biq}) admits the involutions
$$
F_1: \quad (X,Y)\mapsto(Y,X), \qquad
F_2: \quad (X,Y)\mapsto(X,Y^\dagger ),
$$
where $y=Y^\dagger$ is the unique solution to the
quadratic equation $\mathcal{B}(X,y)=0$ apart from $y=Y$. Clearly the
composition $G=F_2\circ F_1$ of these two involutions yields
the mapping of the plane (\ref{plamap}) defined by (\ref{bimap}),
and it is easy to see that the quantity (\ref{kint}) is conserved
by this mapping, so that  $\tilde{K}=K=\mathrm{constant}$ and the points
$(u_{n-1},u_n)$ lie on the level set (\ref{biq})
for all $n$ (see e.g. \cite{bastien}
for more details of this construction).
$G$  is thus a birational map $\hat{E}\to\hat{E}$ composed of
two involutions $F_j$, $j=1,2$. We claim that for
all $P\in\hat{E}$ each $F_j:\,
P\mapsto 2C_j-P$ for some point $C_j\in\hat{E}$, for $j=1,2$
respectively, and hence $G$ is a translation by the 
point $\hat{P}=2C_2-2C_1$, so   
$G:\,P\mapsto P+\hat{P}$ as is required
for  (\ref{pseq}) to hold.

To see why the claim is true, note that each 
involution $F_j$, $j=1,2$ is a 
birational map $\hat{E}\to\hat{E}$ 
with four fixed points. Consider first $F_1$, noting that 
any morphism of 
elliptic curves is a combination of a translation by a point and an isogeny 
$\phi$ 
(see \cite{silver1}, p.75). Since $F_1$ is birational we require  
that $\phi\in\mathrm{Aut}(\hat{E})$, and because it is an involution the only 
possibilities are $\phi =[\pm 1]$. If $\phi =[1]$ then 
we must have $F_1:\, P\mapsto P+\Omega$ where $\Omega$ is a non-identity 
element of the 2-torsion subgroup (of order four),  
but then there are no fixed points. Thus it must be  the case 
that $\phi=[-1]$ and so $F_1:\, P\mapsto 2C_1-P$ where the four possible 
points $C_1+\Omega$ (with $\Omega$ being 2-torsion) 
account for the fixed points, which for 
$F_1$ correspond 
to the roots of $\mathcal{B}(X,X)=0$. The same argument applies to $F_2$. 

To obtain the expression (\ref{efn}), 
observe that because $\hat{E}$ has genus one it can be uniformized using 
elliptic functions; in \cite{ww} this result is attributed to 
Clebsch. Thus we can write $(X,Y)=(g(z),f(z))$ where $g,f$ are 
rational functions of the Weierstrass functions $\wp(z)$, $\wp'(z)$ associated 
with a curve $E$ in Weierstrass canonical form (\ref{canon}) 
with invariants $g_2$, $g_3$, and $z\in\mathrm{Jac}(E)$. 
Now the map defined by 
(\ref{bimap}) corresponds to addition of a point, which is 
equivalent to  
a shift $z\mapsto z+\ka$ on the Jacobian. Then we see that 
for some $z_n$ we have $(u_{n-1},u_n)=(\,g(z_n)\, ,\, f(z_n)\,)$ 
and $(u_n,u_{n+1})=(\,g(z_{n+1})\, ,\, f(z_{n+1})\,) 
=(\,g(z_{n}+\ka  )\, ,\, f(z_{n}+\ka )\,)$, 
whence it follows that $g(z)=f(z-\ka )$ and 
$u_n =g(z_0+(n+1)\ka)=f(z_0+n\ka )$ as required.  
\end{proof} 

\begin{remark} In essence, the above result for 
symmetric biquadratic curves was already known to 
Euler (see Theorem 9 in \cite{what}). A thorough treatment of 
general 
integrable maps associated with elliptic curves over an 
arbitrary field can be found in the recent work \cite{jrv}.  
For our purposes here, we need to construct the explicit uniformization 
for particular families of curves of the above type, and the exact 
form of the corresponding 
expressions in terms of Weierstsrass functions 
changes when one or another of the parameters in (\ref{biq}) 
degenerates to zero. 
For example, we see that 
for the mapping (\ref{nonli}) associated with 
Somos 4 sequences, the invariant
curve 
defined by (\ref{jint})  in 
Proposition 2.1 corresponds to the degenerate case of (\ref{biq}) obtained 
by setting $e=1$, $d=c=0$.  In the proof  of the analogous result for 
Somos 5 sequences, we shall see that the relevant mapping corresponds to 
the case when (\ref{biq}) degenerates to a cubic curve by 
setting $e=c=0$, cf. (\ref{njint}) and  (\ref{necurv}) below.      
\end{remark} 

\begin{theorem} 
Given a Somos 5 sequence satisfying (\ref{s5}), the 
corresponding quantity $h_n$ defined by (\ref{hdef}) 
satisfies the second order nonlinear mapping 
\beq 
\label{hmap} 
h_{n-1}h_nh_{n+1}=\tal h_n +\tbe , 
\eeq 
which has the 
conserved quantity 
\beq \label{njint}   
\tJ = h_{n-1}+h_n +\tal \left( \frac{1}{h_{n-1}} + 
\frac{1}{h_n}\right) +\frac{\tbe}{h_{n-1}h_n}. 
\eeq 
The general solution of this map can be written as 
\beq \label{hsol} 
h_n= \frac{\si (z_0+(n+2)\ka )\,\si (z_0 +(n-1)\ka )} 
{\si (\ka )^4\,\si (z_0+n\ka )\,\si (z_0+(n+1)\ka )} 
=-\frac{\wp '(\ka )}{2} \left( \frac{\wp '(z_0+n\ka ) -\wp' (\ka )} 
{\wp (z_0 +n\ka )-\wp (\ka )} \right) + \frac{\wp ''(\ka )}{2} 
,   
\eeq  
where $\si (z) =\si (z;g_2,g_3)$ and 
$\wp (z) =\wp (z;g_2,g_3)$ are Weierstrass functions 
associated with an elliptic curve 
\beq \label{ecurv} 
E: \qquad y^2=4x^3-g_2x-g_3.  
\eeq  
The invariants of the curve  
are determined from the parameters $\tal$, $\tbe$ together with the 
initial data $h_0$, $h_1$ for the map (\ref{hmap}) according to the 
formulae 
\beq 
\label{lmform} 
\tmu = (\tbe +\tal \tJ )^{1/4},  
\qquad \tla =\frac{1}{3\tmu^2}\left( \frac{\tJ^2}{4}+\tal \right),   
\eeq  
\beq 
\label{gform} 
g_2 =12\tla^2 -2\tJ, \qquad 
g_3 = 4\tla^3 -g_2\tla -\tmu^2.  
\eeq 
The parameters $z_0$, $\ka\in \mathbb{C}$
are then determined 
by  the elliptic integrals 
\beq 
z_0=\pm\int_\infty^{x_0} \frac{\rd x}{y},
\qquad
\ka =\pm\int_\infty^{\tla } \frac{\rd x}{y}.
\label{eints2}
\eeq 
on the curve $E$, 
where 
\beq 
\label{x0} 
x_0=\tla +\frac{\tmu^2}{h_{-1}+h_0-\tJ}, \qquad 
h_{-1}=\frac{\tal h_0 +\tbe}{h_0h_1}  
\eeq 
and the relative sign of $\ka$ and $z_0$ is fixed by the 
requiring the consistency of  
\beq \label{consis} 
\tmu =\wp '(\ka), \qquad 
\wp'(\ka )  \wp '(z_0)= 
(x_0-\tla )(h_{-1}-h_0).
\eeq  
\end{theorem} 


\begin{proof} 
Letting $\tJ_n$ denote the right hand side of 
(\ref{njint}), it is clear that 

\beq \label{jdiff} 
\tJ_{n+1} 
-\tJ_n = \frac{(h_{n+1}-h_{n-1})}{h_{n-1}h_nh_{n-1}} 
(h_{n+1}h_nh_{n-1}-\tal h_n-\tbe ),  
\eeq 
so assuming that $h_{n+1}\neq h_{n-1}$ for all $n$ then 
$\tJ$ is conserved if and only if $h_n$ satisfies 
(\ref{hmap}). Alternatively note that rewriting  
(\ref{j3int}) in terms of $h_n$ gives precisely 
(\ref{njint}), and rearranging 
the latter expression for $\tJ$ yields 
\beq \label{hcurv} 
(h_{n-1}+h_n-\tJ )(h_{n-1}h_n +\tal )+\tbe +\tal \tJ =0, 
\eeq 
which implies that for any sequence $h_n$, $n\in\mathbb{Z}$  
satisfying (\ref{hmap}), 
$(X,Y)=(h_{n-1},h_n)$ is a point on the curve 
\beq \label{necurv} 
(X+Y-\tJ )(XY+\tal )+\tbe +\tal \tJ =0 
\eeq 
for all $n$.

We claim that this elliptic 
curve is birationally equivalent to (\ref{ecurv})      
via the transformation 
\beq \label{birat} 
X=\frac{\tmu}{2}\left(\frac{y+\tmu }{x-\tla}\right)+\frac{\tJ}{2}, 
\qquad 
Y=-\frac{\tmu}{2}\left(\frac{y-\tmu }{x-\tla}\right)+\frac{\tJ}{2},   
\eeq  
where the parameters $\tla$, $\tmu$ and invariants $g_2$, $g_3$ 
are found from the coefficients of the curve (\ref{necurv}) 
using the formulae (\ref{lmform}) and (\ref{gform}). 
To verify this, it is instructive to use the analytic formulae (\ref{hsol}) 
which provide a uniformization of (\ref{necurv}). First note that      
the two formulae for $h_n$ in (\ref{hsol}) are easily seen to be 
equivalent by noting that both expressions are elliptic functions of 
$z_0$ with simple zeros at $z_0=-(n+2)\ka$, $-(n-1)\ka$ 
and simple poles at $z_0=-n\ka$, $-(n+1)\ka$, with the same residues at 
these poles. Similarly, shifting $n\to n-1$ in the first (sigma function) 
expression in (\ref{hsol})  
gives 
\beq \label{hmsol}
h_{n-1}= \frac{\si (z_0+(n+1)\ka )\,\si (z_0 +(n-2)\ka )}
{\si (\ka )^4\,\si (z_0+(n-1)\ka )\,\si (z_0+n\ka )}
=\frac{\wp '(\ka )}{2} \left( \frac{\wp '(z_0+n\ka ) +\wp' (\ka )} 
{\wp (z_0 +n\ka )-\wp (\ka )} \right) + \frac{\wp''(\ka )}{2} 
, 
\eeq
where once again these two formulae for $h_{n-1}$ are seen to be 
equivalent by looking at their poles and zeros. Then we can make the 
identifications 
\beq \label{idparam} 
\tla =\wp(\ka ), \quad \tmu =\wp'(\ka )=(\tbe +\tal\tJ )^{1/4}, 
\quad \tJ =\wp''(\ka ), 
\quad \tal =-\wp'(\ka )^2\Big(\wp (2\ka )-\wp (\ka )\Big),    
\eeq  
which are consistent with (\ref{lmform}) and (\ref{gform}), i.e.   
$(\tla ,\tmu)=(\wp (\ka ),\wp'(\ka ))$ is a point on the curve 
(\ref{ecurv}).   

Taking the sum of the right hand ($\wp$ function) expressions 
in (\ref{hsol}) and (\ref{hmsol}) 
gives 
\beq \label{sum} 
h_{n-1}+h_n =\wp '' (\ka )  
+\frac {\wp'(\ka )^2}{\wp (z_0+n\ka )-\wp (\ka )}, 
\eeq 
which by $(\ref{idparam})$ is equivalent to 
$X+Y=\tJ +\tmu^2/(x-\tla )$ with the identification $(x,y)= 
(\wp (z_0 +n\ka ),\wp'(z_0 +n\ka ))$. On the other hand,  
taking the product of the left hand (sigma function) 
expressions
in (\ref{hsol}) and (\ref{hmsol})
gives
\beq \label{prod}
h_{n-1}h_n =\frac {\si (z_0+(n+2)\ka )\,\si (z_0 +(n-2)\ka )}
{\si (\ka )^8\,\si (z_0+n\ka )^2} 
=\wp'(\ka )^2\Big(\wp (2\ka )-\wp (z_0+n\ka )\Big),
\eeq 
where we have used (\ref{addweier}) and  
the first identity in (\ref{efns}). 
With the same identifications as before, 
this  product identity implies that we have 
$XY+\tal = -\tmu^2 (x-\tla )$, 
and hence $(X+Y-\tJ )(XY+\tal )=-\tmu^4=-(\tbe +\tal\tJ )$ which 
means that $(X,Y)=(h_{n-1},h_n)$ lies on the curve 
(\ref{necurv})  for any $n$, as claimed. In turn this means that 
the sequence $\{ h_n \,|\, n\in\mathbb{Z}\}$ defined by 
(\ref{hsol}) has $\tJ =\wp ''(\ka )$ as a conserved quantity, and hence 
by (\ref{jdiff}) it must satisfy the second order map (\ref{hmap}). 

To see that (\ref{hsol}) represents the unique solution 
of the initial value problem for the map (\ref{hmap}), observe that 
it depends on the four parameters $g_2,g_3,z_0,\ka$. Up to 
the involution $z_0\to -z_0$, $\ka \to -\ka$ (which does not change 
the solution), all of  the latter are   
uniquely determined from the initial data $h_0$, 
$h_1$ and parameters $\tal$, $\tbe$ by first constructing the  
first integral $\tJ =\tJ_0$ (i.e. setting $n=0$ in (\ref{njint})) 
and then using the formulae (\ref{lmform}) and (\ref{gform}) to 
find the invariants. Finally the base point $z_0$ and the shift 
$\ka$ on the Jacobian Jac($E$) of the curve (\ref{ecurv}) 
are obtained by 
inverting $x_0=\wp(z_0)$, $\tla =\wp(\ka )$ to yield the elliptic 
integrals (\ref{eints2}), where $x_0$ in (\ref{x0}) 
is found by rearranging (\ref{sum}) 
for $n=0$. The consistency requirements (\ref{consis}) 
follow from fixing the branch of fourth root taken for $\tmu$ 
in (\ref{lmform}) and ensuring that the value of $\wp'(z_0 )$ 
agrees with this in the formulae for $h_0$ and $h_{-1}$ (that is, 
the $\wp$ function expressions in (\ref{hsol}) and (\ref{hmsol}) 
when $n=0$). 
\end{proof}

\begin{prop} 
For any solution $\tau_n$, $n\in\mathbb{Z}$ 
of a Somos 5 recurrence (\ref{s5}), the subsequences $\tau_n^* = 
\tau_{2n}$ and $\tau_n^* = 
\tau_{2n+1}$, of respectively even and odd index terms, both satisfy 
the Somos 4 recurrence 
\beq
\tau_{n+2}^*\, 
\tau_{n-2}^* 
=\al^* \,
\tau_{n+1}^* 
\tau_{n-1}^* 
+\beta^* 
\,
(\tau_n^*)^2,
\label{starbil}
\eeq
where 
\beq 
\label{abstar} 
\al^* = \tbe^2, 
\quad 
\be^* = \tal (2\tbe^2+\tal\tbe\tJ +\tal^3),  
\eeq  
with $\tJ$ being the associated integral (\ref{njint}) of the 
map (\ref{hmap}). 
\end{prop}

\begin{proof} 
Looking at the formula (\ref{prod}) we see that 
$$ 
h_nh_{n-1}=f_{n+1}(f_n)^2f_{n-1}=\frac{\tau_{n+2}\tau_{n-2}}{(\tau_n)^2} 
=\wp '(\ka )^2\Big( \wp (2\ka )-\wp (z_0+n\ka )\Big),   
$$ 
which is almost identical to the formula (\ref{desol}) for the solution 
of the map (\ref{nonli}), except for the scale factor 
$\tmu^2 =\wp '(\ka )^2 =\sqrt{\tbe +\tal\tJ }$. Using the 
scaling property  
$$ 
\wp (z;g_2,g_3)=\tmu^{-2}\wp (z/\tmu ;\tmu^4g_2,\tmu^6g_3), 
\quad 
\si (z;g_2,g_3)=\tmu \si (z/\tmu ;\tmu^4g_2,\tmu^6g_3), 
$$ 
of the Weierstrass functions, we can rewrite everything 
in terms of rescaled sigma and $\wp$ functions with 
invariants $g_2^*=\tmu^4g_2$, 
$g_3^*=\tmu^6g_3$,   setting $v=\ka/\tmu$ and 
$u_0=z_0/\tmu$ 
to obtain the canonical form 
\beq \label{cansol} 
\frac{\tau_{n+2}\tau_{n-2}}{(\tau_n)^2}=\wp (2v )-\wp (u_0 +nv) 
=\frac{\si (u_0 +(n+2)v) \si (u_0 +(n-2)v)}{\si (u_0+nv)^2\si (2v)^2}.   
\eeq 
It follows immediately from Proposition 2.1 that both 
$F_n^*=\tau_{2n+2}\tau_{2n-2}/(\tau_{2n})^2$ 
and 
$F_n^*=\tau_{2n+3}\tau_{2n-1}/(\tau_{2n+1})^2$ are solutions of the nonlinear  
map 
$$ 
F_{n+1}(F_n)^2F_{n-1}=\al^*F_n+\be^*, 
$$ 
with 
\beq \label{dup}  
\al^*=\wp '(2v)^2, \qquad 
\be^* = \wp '(2v)^2 \Big(  
\wp (4v) -\wp (2v)   
\Big) .   
\eeq 
The required result follows immediately, upon noting that 
the duplication formula for the Weierstrass $\wp$ function 
can be used in (\ref{dup}) 
to 
rewrite $\wp (4v)$ in terms of $\wp (2v)$ and its derivatives, 
and similarly for $\wp (2v)$ in terms of $\wp (v)$. After  
further tedious computations, scaling back to the 
original Weierstrass functions $\tal =\wp (\ka )$, 
$\tmu = \wp'(\ka )$, $\tJ=\wp''(\ka )$ with invariants $g_2$, $g_3$, 
we arrive at the expressions  
(\ref{abstar}). 
\end{proof}

\begin{theorem}    
The general solution of the Somos 5 recurrence takes the alternating 
form  
\beq 
\label{final} 
\tau_{2n}=A_+ \,B_+^n \frac{\si (u_0 +2nv)}{\si (2v)^{n^2}}, 
\qquad 
\tau_{2n+1}=A_- \,B_-^n \frac{\si (u_0 +(2n+1)v)}{\si (2v)^{n^2}}, 
\eeq 
where here $\si (u)=\si (u;g_2^*,g_3^*)$ 
denotes the Weierstrass sigma function associated 
with  the elliptic curve 
\beq 
E^*: \qquad y^2=4x^3-g_2^*x-g_3^*, 
\label{estar} 
\eeq 
where $u_0$, $v$ and the invariants are given in terms of the parameters 
(\ref{lmform}), (\ref{gform}) in Theorem 2.7 by 
$$ 
g_2^*=\tmu^4 \,g_2, \qquad  
g_3^*=\tmu^6 \,g_3, \qquad  
u_0 =z_0/\tmu, \qquad  
v =\ka /\tmu,    
$$  
and the prefactors are   
\beq  
\label{apm}  
A_+=\frac{\tau_0}{\si (u_0)}  ,\qquad A_-=  
\frac{\tau_1}{\si (u_0+v)}  \eeq  
\beq 
B_+=\frac{\si (2v)\si (u_0)\tau_2}{\si (u_0+2v)\tau_0}  , 
\qquad B_-=\frac{\si (2v)\si (u_0+v)\tau_3}{\si (u_0+3v)\tau_1}  . 
\label{bpm} 
\eeq 
The solution (\ref{final}) depends on seven independent parameters 
$g_2^*,g_3^*,u_0,v,A_+,A_-$ and $B_+$, with $B_{\pm}$ defined  
in (\ref{bpm}) being related by 
\beq \label{brel} 
\frac{B_+}{B_-}=-\si (2v). 
\eeq 
\end{theorem} 


\begin{proof} 
Starting from the result of Proposition 2.8, 
the formula (\ref{cansol}) for the ratios of the $\tau_n$ can be 
used to solve for the even and odd index subsequences separately, since 
these are both Somos 4 sequences. In that case, Theorem 1.1 applies   
separately to the two  subsequences 
$\tau_{2n}$ and $\tau_{2n+1}$, which yields (\ref{final}) -  
the proof for each subsequence is the same as in \cite{honeblms}.       
However, the solution can depend on at most 7 parameters, as is 
confirmed by  the 
relation (\ref{brel}) between $B_\pm$. To prove the 
latter relationship, it suffices to rewrite the ratios of sigmas 
and $\tau_n$ in terms of $h_1$, and see that most terms cancel, 
taking care to scale appropriately with 
\beq \label{tmu} 
\tmu =\wp '(\ka ;g_2,g_3), 
\qquad 
\tmu^4=\wp '(v;g_2^*,g_3^*). 
\eeq 
\end{proof} 

\begin{remark} 
The rescaled curve (\ref{estar}) is  
more convenient when working in fields other than $\mathbb{C}$,  
since then one does not need to define the curve over the field  
extension by the fourth root $\tmu = (\tbe +\tal \tJ )^{1/4}$.   
\end{remark}

\begin{cor} 
The solution of the third order 
nonlinear map (\ref{nonli}) can be written in the 
alternating form 
\beq \label{n3sol} 
f_{2n}=f_0 \left(\frac{\wp (v)-\wp(u_0+2nv)}{\wp (v)-\wp(u_0)}\right), 
\qquad 
f_{2n+1}=f_1 \left(\frac{\wp (v)-\wp(u_0+(2n+1)v)}{\wp (v)-\wp(u_0+v)}\right), 
\qquad 
\eeq 
where $\wp(u)=\wp(u;g_2^*,g_3^*)$ is the Weierstrass $\wp$ function 
associated with the curve (\ref{estar}), with the parameters $u_0$, $v$  
as given in Theorem 2.9.  
\end{cor}

\begin{proof} 
This follows by substituting 
(\ref{final}) into (\ref{fdef}) and rewriting suitable combinations of 
the prefactors 
$A_\pm$, $B_\pm$ in terms of $f_0$ and $f_1$.  
\end{proof}

\begin{cor} 
The    
terms $\tau_n$, $n\in\mathbb{Z}$  
of a Somos 5 sequence satisfy the Hankel 
determinant formula  
\beq 
\label{hankel4}
a_1 a_2 \tau_{n+m+1}\tau_{n-m}=
\left|\begin{array}{cc}
a_{m+1} \tau_{n+2} & a_{m-1} \tau_n \\
a_{m+2}\tau_{n+1} & a_m \tau_{n-1}
\end{array} \right|
\eeq 
for all $m,n\in \mathbb{Z}$, where the quantities $a_m$, 
which form an associated elliptic divisibility sequence,  
are given by 
$$a_m=\frac{\si (m\ka )}{\si (\ka )^{m^2}}  
, \qquad 
m\in \mathbb{Z}$$ in terms of the sigma function $\si (z)=\si (z;g_2,g_3)$  
for the  
curve (\ref{ecurv}) in Theorem 2.7.  
\end{cor}

\begin{proof} 
Using the scaling by $\tmu$ as in (\ref{tmu}), 
the formulae (\ref{final}) can be rewritten as  
\beq 
\label{final2} 
\tau_{2n}=A_+ \,B_+^n \frac{\tmu^{n^2-1}\si (z_0 +2n\ka )}{\si (2\ka )^{n^2}}, 
\qquad 
\tau_{2n+1}=A_- \,B_-^n \frac{\tmu^{n^2-1} 
\si (z_0 +(2n+1)\ka )}{\si (2\ka )^{n^2}}, 
\eeq 
in terms of the sigma function for the original curve (\ref{ecurv}) as 
in Theorem 2.7. In terms of this sigma function   
it is easy to show, upon rescaling (\ref{brel}) by 
$\tmu =\wp'(\ka ;g_2,g_3 ) 
=-\si (2\ka ;g_2,g_3)/\si (\ka ;g_2,g_3)^4$ (compare with the 
first formula of (\ref{efns})), that 
$$ 
B_+/B_-=\si (\ka ;g_2,g_3)^4. 
$$ 
For verification of    
(\ref{hankel4}), the symmetry $m\to -m-1$ means 
that it is sufficient to check only the   
cases when $m,n$ are either both odd or both even, to take the 
different expressions (\ref{final2}) 
for $\tau_n$ with odd/even $n$ into account. Thereafter the proof 
is the same as for Corollary 1.3.  
\end{proof}

\begin{cor}  
The terms of a Somos 5 sequence, 
satisfying the recurrence (\ref{s5}), have the 
leading order asymptotic behaviour 
\beq \label{asymp} 
\log |\tau_n| \sim \left(\mathrm{Re}\, \left\{\frac 
{\eta_1 v^2}{2\om_1}\right\} -\frac{1}{4}\log |\si (2v)|\right)n^2 
,\qquad n\to \infty  
\eeq 
for $v\in\omega_1\mathbb{R}$, 
where $\omega_1$ is a half-period and  $\eta_1=\zeta (\om_1)$, with  
the invariants $g_2^*$, $g_3^*$ being as in Theorem 2.9. 
\end{cor}

\begin{proof} 
The sigma function is given by the product formula 
$$ 
\si (z) = \frac{2\om_1}{\pi}\exp \left(\frac{\eta_1z^2}{2\om_1}\right) 
\,\sin \left( \frac{\pi z}{2\om_1} \right) \, 
\prod_{n=1}^\infty 
\left\{\frac{1-2q^{2n}\cos(\pi z/\om_1)+q^{4n}}{(1-q^{2n})^2}\right\} 
$$ 
with $q=\exp (\pi i\om_2/\om_1 )$  (see $\S$20.421 and also 
$\S$21.43 in \cite{ww} for the corresponding theta function expression). 
The result then follows  
upon taking the logarithm of (\ref{final}). The corrections 
to (\ref{asymp}) at $O(n)$ can be calculated similarly, although these 
have a different form for even/odd $n$. 
\end{proof}  

\section{Example: Somos (5)}

To illustrate the above results, we consider the example of  
the Somos (5) sequence (\ref{somos5}). The initial data and 
parameters as in (\ref{s5data}) give $h_0=2$, $h_1=1$ for the  
second order map 
(\ref{hmap}), whence $h_{-1}=3/2$ with 
$\tal =1=\tbe$, 
so applying the formulae in Theorem 2.5 we have   
$$ 
\tJ =5, \quad \tmu =6^{1/4}, \quad \tla =\frac{29}{12\sqrt{6}}, 
\quad x_0 =-\frac{19}{12\sqrt{6}},  
$$ 
whence 
$g_2=121/72$, $g_3=-845/(1296\sqrt{6})$. 
Rescaling by appropriate factors of $\tmu$, as in Theorem 2.9, 
we arrive at the 
curve 
\beq \label{s5e} 
E^*: \qquad y^2=4x^3-\frac{121}{12}x+\frac{845}{216} 
, \qquad g_2^*=\frac{121}{12}, \quad g_3^*=-\frac{845}{216}, 
\eeq 
which has the $j$-invariant 
$$ 
j=\frac{1728(g_2^*)^3}{(g_2^*)^3-27(g_3^*)^2}=\frac{1771561}{612}. 
$$ 
Hence the Somos (5) sequence (\ref{somos5}) corresponds to the 
sequence of points 
\beq \label{points} 
(x_0^*,y_0^*)+n(\la^* ,\mu^*)=(-19/12,2)+n(29/12,6)= 
(5/12,0)+(n-2)\,(29/12,6) 
\eeq 
on the rescaled curve (\ref{s5e}). 
   
For comparison with results of other authors, we note that van der Poorten 
\cite{vdp} has obtained the Somos (5) sequence from the minimal model 
$V^2+UV+6V=U^3+7U^2+12U$, while  Zagier \cite{zagier} 
found the related curve 
$y^2+xy=x(x-1)(x+2)$; both of these curves have the same $j$-invariant 
as (\ref{s5e}), and hence they are birationally equivalent to each 
other. Zagier's discussion makes use of the inspired substitution 
(\ref{hdef}) for the Somos (5) sequence (\ref{somos5}), although 
the presentation of the associated curve in 
the form (\ref{necurv}) contains 
a typographical error in \cite{zagier}. Other approaches  
to the sequence (\ref{somos5}) can be found on postings to Propp's 
``robbins'' forum 
\cite{propp} - see in particular 
the references to unpublished work of Elkies 
in \cite{heron}, where the sequence (\ref{somos5}) is related 
to the problem of Heron triangles with two rational medians.  

With the curve (\ref{s5e}) and the sequence of points (\ref{points}) 
we can use Theorem 2.9 to calculate the general term of the 
Somos 5 recurrence in the form (\ref{final}) 
with 
$$ 
u_0 =0.163392411+\om_2, \qquad 
v=-0.672679183, 
$$ 
where the real and imaginary half-periods are 
$$ 
\om_1=1.181965956, 
\qquad 
\om_2=0.973928783\, i 
$$ 
respectively (doing all the elliptic integrals to 9 decimal 
places with MAPLE). The result of Corollary 2.13 then implies that 
the terms of the sequence (\ref{somos5}) have the leading order 
asymptotics 
$$ 
\log \tau_n \sim  
0.071626946 \, n^2, \qquad n\to\infty . 
$$ 

\section{Conclusions} 

The theorems presented in \cite{honeblms} and above 
provide a simple algorithmic approach to 
determining elliptic curves associated with Somos 4 and Somos 5 sequences, 
and finding the complete solution of the associated initial value problems. 
Our approach is based on making the connection with explicitly solvable maps 
with conserved quantities  
i.e. (\ref{nonli}) and (\ref{hmap}). A symplectic 
map with sufficiently many conserved quantities 
is called integrable 
\cite{rag, ves}. As they stand the maps presented here are not symplectic, 
but by making a change of variables, 
a suitable symplectic structure was given 
in \cite{honeblms}. 

The results presented above naturally lead to the question of 
how to solve Somos $(2N+2)$ recurrences with $N+2$ terms, of the form 
$$ 
\tau_{n+N+1}\tau_{n-N-1}=\sum_{j=0}^N \al_j \,\tau_{n+j}\tau_{n-j}, 
$$ 
with constant coefficients $\al_j$, and similarly for Somos $(2N+3)$. 
By using the continued fraction expansion of the square root of 
a sextic, van der Poorten \cite{vdp2} 
has generated  a 
certain class of Somos 6 sequences associated 
with  genus 2 curves. With a different approach \cite{beh} based on 
the addition formulae for Kleinian sigma functions (see \cite{bel} 
and references), we have constructed Somos 8 sequences derived from  
genus 2 curves, taking a different (quintic) 
affine model compared with van der 
Poorten. The recurrences in \cite{beh} generalise  
the genus 2 case of  
Cantor's hyperelliptic division polynomials \cite{cantor} 
(for analytic formulae, see \cite{matsupsi}), and 
they are related to a family of integrable symplectic maps 
(discrete H\'{e}non-Heiles systems). We should also remark that 
Buchstaber and Krichever have derived bilinear addition formulae for 
Riemann theta functions \cite{bk}, 
which have exactly $N+2$ terms in genus $N$. 
Whether this result can lead to an effective solution to the initial 
value problem for these higher order recurrences is the subject of 
further investigation. 
It would also be interesting to see
how some of our results on Somos 4 and 5 sequences could be modified
in the $p$-adic setting \cite{silverman}.

\noindent {\bf Acknowledgments.}
I am grateful to Graham Everest for introducing me to the  
arithmetic of quadratic recurrence sequences, to  
Christine Swart for sending me her thesis \cite{swart}, 
and to Marc Rogalski for sending me his work \cite{bastien}.   
Thanks also to Tom Ward,  
David Cantor and Alf van der Poorten for interesting  
correspondence on related matters,   
and to the University of Kent  
for supporting  
the project  
{\it Algebraic curves and functional equations 
in mathematical physics} with a Colyer-Fergusson Award.   
Finally, I would like to thank the referee for suggesting  
Proposition 2.5.

\bibliographystyle{amsplain}

\end{document}